\documentclass{amsart}
\usepackage{amsmath}
\usepackage{graphicx}
\usepackage{amsfonts}
\usepackage{amssymb}
\newtheorem{theorem}{Theorem}
\theoremstyle{plain}

\newtheorem{corollary}{Corollary}

\newtheorem{definition}{Definition}
\newtheorem{example}{Example}

\newtheorem{remark}{Remark}

\numberwithin{equation}{section}

\newcommand{\abs}[1]{\left\lvert#1\right\rvert}
\DeclareMathOperator{\supp}{supp}

\DeclareMathOperator{\conv}{conv}
\DeclareMathOperator{\clco}{\overline{conv}}
\DeclareMathOperator{\cone}{cone}
\DeclareMathOperator{\clcone}{\overline{cone}}

\begin{document}
\title{The proof of Tchakaloff's Theorem}
\author{Christian Bayer, Josef Teichmann}
\address{Technical University of Vienna, e105, Wiedner Hauptstrasse 8-10, A-1040 Wien, Austria}
\email{cbayer@fam.tuwien.ac.at, jteichma@fam.tuwien.ac.at}
\thanks{The authors are grateful to Prof.~Peter Gruber for mentioning the word
\textquotedblright St\"{u}tzebene\textquotedblright\ in the right moment. The
first author acknowledges the support from FWF-Wissenschaftskolleg
\textquotedblright Differential Equations\textquotedblright\ W 800-N05. The
second author acknowledges the support from the RTN network HPRN-CT-2002-00281
and from the FWF grant Z-36.}

\begin{abstract}
We provide a simple proof of Tchakaloff's Theorem on the existence of
cubature formulas of degree $m$ for Borel measures with moments up to
order $m$. The result improves known results for non-compact support,
since we do not need conditions on $(m+1)$st moments. In fact we reduce
the classical assertion of Tchakaloff's Theorem to a well-known
statement going back to F. Riesz.
\end{abstract}\maketitle

We consider the question of existence of cubature formulas of degree $m$ for
\emph{Borel measures} $\mu$, i.~e.~a measure defined on the Borel $\sigma$-algebra,
where moments up to degree $m$ exist:

\begin{definition}
Let $\mu$ be a positive Borel measure on $\mathbb{R}^{N}$ and $m\geq1$ such
that
\[
\int_{\mathbb{R}^{N}}\left\|  x\right\|  ^{k}\mu(dx)<\infty
\]
for $0\leq k\leq m$ holds true. A \emph{cubature formula of degree }$m$ for
$\mu$ is given by an integer $k\geq1$, points $x_{1},\dots,x_{k}\in\supp\mu$,
weights $\lambda_{1},\dots,\lambda_{k}>0$ such that
\[
\int_{\mathbb{R}^{N}}P(x)\mu(dx)=\sum_{i=1}^{k}\lambda_{i}P(x_{i})
\]
for all polynomials on $\mathbb{R}^{N}$ with degree less or equal $m$, where
$\supp \mu$ denotes the closed support of the measure $\mu$, i.~e.~the complement
of the biggest open set $O \subset \mathbb{R}^N$ with $\mu(O)=0$.
\end{definition}

Cubature formulas of degree $m$ have been proved to exist for Borel measures
$\mu$, where the $(m+1)$st moments exist, see \cite{cur/fia:02} and
\cite{put:97}. The result in the case of compact $\supp\mu$ is classical, and
due to Tchakaloff (see \cite{tch:57}), hence we refer to the assertion as
Tchakaloff's Theorem.

We collect some basic notions and results from convex analysis, see for
instance \cite{roc:72}: fix $N\geq1$, for some set $S\subset\mathbb{R}^{N}$
the convex hull of $S$, i.~e.~the smallest convex set in $\mathbb{R}^{N}$
containing $S$, is denoted by $\conv(S)$, the (topological) closure of $\conv
(S)$ is denoted by $\clco(S)$. The closure of a convex set is convex. Note
that the convex hull of a compact set is always closed, but there are closed
sets whose convex hull is no longer closed (see \cite{roc:72}).

Closed convex sets can also be described by their \emph{supporting
hyperplanes}. Given a convex set $C$. Let $y\in\partial C:=\overline
{C}\setminus\operatorname*{int}(C)$ be a boundary point. There is a linear
functional $\mathit{l}_{y}$ and a real number $\beta_{y}$ such that the
hyperplane defined by $\mathit{l}_{y}=\beta_{y}$ contains $y$, and $C$ is
contained in the closed half-space $\mathit{l}_{y}\leq\beta_{y}$. Hyperplanes
and half-spaces with this property are called supporting hyperplanes and
supporting half-spaces, respectively. Moreover, $\overline{C}$ is the
intersection of all its supporting half-spaces. Furthermore, if $C$ is not
contained in any hyperplane of $\mathbb{R}^{N}$ (i.e. it has non-empty
interior), then a point $x\in C$ is contained in a supporting hyperplane if
and only if $x\in C$ is not an interior point of $C$ (see \cite{roc:72}, Th.
11.6). This means that we can characterize the boundary of $C$ as those
points, which lie at least in one supporting hyperplane of $C$.

In the case of a convex cone $C$ the supporting hyperplanes can be chosen to
be homogeneous, i.~e.~to be of the form $\mathit{l}_{y}=0$. We denote the
convex cone generated by some set $A\subset\mathbb{R}^{N}$ by $\cone(A)$ and
its closure by $\clcone(A)$.

We also introduce the notion of the \emph{relative interior} of a convex set
$C$: a point $x\in C$ lies in the relative interior $\operatorname*{ri}(C)$ if
for every $y\in C$ there is $\epsilon>0$ such that $x-\epsilon(y-x)\in C$. In
particular we have that the relative interior of a convex set $C$ coincides
with the relative interior of its closure $\overline{C}$. Interior points of
$C$ lie in the relative interior (see \cite{roc:72}), this remains true even
if $C$ lies in an affine subspace of $\mathbb{R}^{N}$, and a point of $C$ lies
in the interior with respect to the subspace topology.

Given a measure $\mu$ on some measurable space $(\Omega,\mathcal{F})$ and a Borel measurable map
$\phi:\Omega\rightarrow\mathbb{R}^{N}$, we denote by $\phi_{\ast}\mu$
the \emph{push-forward Borel measure} on $\mathbb{R}^{N}$, which is defined via
\[
\phi_{\ast}\mu(A):=\mu(\phi^{-1}(A)),
\]
for all Borel sets $A\subset\mathbb{R}^{N}$.

\begin{theorem}
Let $\mu$ be a positive Borel measure on $\mathbb{R}^{N}$, such that the first
moments exist, i.~e.
\[
\int_{\mathbb{R}^{N}}\left\|  x\right\|  \mu(dx)<\infty,
\]
and let $A\subset\mathbb{R}^{N}$ be a measurable set with $\mu(\mathbb{R}%
^{N}\setminus A)=0$. Then the first moment $E=\int_{\mathbb{R}^{N}}x\mu(dx)$,
where $x$ denotes the vector $(x_{1},\dots,x_{N})$, is contained in
$\cone(A)$.\label{cubature result}
\end{theorem}

\begin{proof}
We first assume that there is no $B\subset A$ with $\mu(A\setminus B)=0$ such
that $B$ is contained in a hyperplane, since otherwise we could work in a
lower-dimensional space instead (with $A$ replaced by $B$). Fix some
$y\in\overline{K}\setminus\operatorname*{int}(K)$ in the boundary of
$K=\cone(A)$. Then all linear functionals $\mathit{l}_{y}:\mathbb{R}%
^{N}\rightarrow\mathbb{R}$ corresponding to the supporting half-spaces
$\mathit{l}_{y}\leq0$ at $y$ are certainly integrable and we have
\[
\mathit{l}_{y}(E)=\int_{\mathbb{R}^{N}}\mathit{l}_{y}(x)\mu(dx)\leq0,
\]
consequently $E\in\clcone(A)$.

By existence of the first moments, for each $\delta>0$ we have $\mu
(\mathbb{R}^{N}\setminus B(0,\delta))<\infty$, where $B(0,\delta)$ denotes the
centered ball with radius $\delta$. Given $\mathit{l}_{y}$ as above, we may
conclude that $\mu(\{x\in A|\mathit{l}_{y}(x)<0\})>0$, since otherwise the
complement in $A$ of the intersection of $A$ with the hyperplane
$\mathit{l}_{y}=0$ would have measure $0$, a contradiction to the assumption
above. Then we can find $\epsilon>0$ such that $0<\mu(\{x\in A|\mathit{l}%
_{y}(x)\leq-\epsilon\})<\infty$ and get
\[
\mathit{l}_{y}(E)=\int_{\mathbb{R}^{N}}\mathit{l}_{y}(x)\mu(dx)\leq
-\epsilon\mu(\{x\in A|\mathit{l}_{y}(x)\leq-\epsilon\})<0.
\]
Hence $E\in\clcone(A)$ is an interior point of $\clcone(A)$. In particular
$E\in\cone(A)$, since the interior lies in the convex cone hull of $A$. If the
first condition is not satisfied, we obtain that $E$ is an interior point of
$\cone(A)$ in an affine subspace of $\mathbb{R}^{N}$ (where the first
condition is satisfied), but then $E$ lies in the relative interior of
$\cone(A)$ in $\mathbb{R}^{N}$, which is the desired result.
\end{proof}

\begin{corollary}
Let $\mu$ be a positive Borel measure on $\mathbb{R}^{N}$ concentrated in
$A\subset\mathbb{R}^{N}$, i.~e.~$\mu(\mathbb{R}^{N}\setminus A)=0$, such that
the first moments exist, i.~e.
\[
\int_{\mathbb{R}^{N}}\left\|  x\right\|  \mu(dx)<\infty.
\]
Then there exist an integer $1\leq k\leq N$, points $x_{1},\dots,x_{k}\in A$
and weights $\lambda_{1},\dots,\lambda_{k}>0$ such that
\[
\int_{\mathbb{R}^{N}}f(x)\mu(dx)=\sum_{i=1}^{k}\lambda_{i}f(x_{i})
\]
for any monomial $f$ on $\mathbb{R}^{N}$ of degree $1$.
\end{corollary}

\begin{proof}
The corollary follows immediately from Theorem~\ref{cubature result} and
Caratheodory's Theorem (see \cite{roc:72}, Th.~17.1 and Cor.~17.1.2).
\end{proof}

\begin{corollary}
Let $\mu$ be a positive measure on the measurable space $(\Omega
,\mathcal{F})$ concentrated in $A\in\mathcal{F}$, i.~e.~$\mu(\Omega\setminus
A)=0$, and $\phi:\Omega\rightarrow\mathbb{R}^{N}$ a Borel measurable map.
Assume that the first moments of $\phi_{\ast}\mu$ exist, i.~e.
\[
\int_{\mathbb{R}^{N}}\left\|  x\right\|  \phi_{\ast}\mu(dx)<\infty.
\]
Then there exist an integer $1\leq k\leq N$, points $\omega_{1},\dots
,\omega_{k}\in A$ and weights $\lambda_{1},\dots,\lambda_{k}>0$ such that
\[
\int_{\Omega}\phi_{j}(\omega)\mu(d\omega)=\sum_{i=1}^{k}\lambda
_{i}\phi_{j}(\omega_{i})
\]
for $1\leq j\leq N$, where $\phi_{j}$ denotes the $j$-th component of $\phi$.
\label{cubature result cooked down}
\end{corollary}

\begin{remark}
In other words, $A\in\mathcal{F}$ such that $\mu(\Omega\setminus A)=0$
correspond to $B\subset\phi(\Omega)$ such that ${\phi}_{\ast}\mu
(\mathbb{R}^{N}\setminus B)=0$.
\end{remark}

\begin{remark}
Note that $\mu(\Omega)=\infty$ is also possible, since we only speak about integrability of $N$ measurable
functions $\phi_{1},\dots,\phi_{N}$. If we have $\mu(\Omega)<\infty$, we could add $\phi_{N+1}=1$, and we obtain in
particular $\sum_{i=1}^{k'}\lambda'_{i}=\mu(\Omega)$ (with possibly different number $ 1 \leq  k' \leq N+1 $ of
points $ x'_i $ and weights $ \lambda'_i $).

In the setting of Theorem~\ref{cubature result} assume that $\mu$ is a probability measure on $\mathbb{R}^{N}$.
Then -- by the previous consideration -- $E = \int_{\mathbb{R}^N}x\mu(dx)$ lies in the convex hull $\conv(A)$. This
fact is well-known in financial mathematics, since it means that the price range of forward contracts is given by
the relative interior of the convex hull of the no-arbitrage bounds of the (discounted) price process (see for
instance~\cite{foe/sch:02}, Th.~1.40).

The result is also well-known in the field of geometry of the moment
problem, see for instance~\cite{kem:65}. As mentioned therein, the
result for compactly supported measures essentially even goes back to
F. Riesz, see~\cite{rie:11}.
\end{remark}

\begin{proof}
We solve the problem with respect to $\phi_{\ast}\mu$ on $\mathbb{R}^{N}$ and
obtain $1\leq k\leq N$, $y_{1},\dots,y_{k}\in\phi(A)$ and $\lambda_{1}%
,\dots,\lambda_{k}>0$ such that
\[
\int_{\mathbb{R}^{N}}f(y)(\phi_{\ast}\mu)(dy)=\sum_{i=1}^{k}\lambda_{i}%
f(y_{i})
\]
for all polynomials $f$ of degree $1$. Thus we obtain points $\omega_{1}%
,\dots,\omega_{k}$ with $\phi(\omega_{i})=y_{i}$ for $1\leq i\leq k$,
furthermore
\[
\int_{\mathbb{R}^{N}}f(y)(\phi_{\ast}\mu)(dy)=\int_{\Omega}(f\circ\phi
)(\omega)\mu(d\omega)
\]
by definition, hence the result.
\end{proof}

In an adequate algebraic framework the previous Theorem \ref{cubature result}
yields all cubature results in full generality, and even generalizes those
results (see \cite{cur/fia:02}, \cite{put:97} and \cite{rez:92} for related
theory and interesting extensive references).

For this purpose we consider polynomials in $N$ (commuting) variables
$e_{1},\dots,e_{N}$ with degree function $\deg(e_{i}):=k_{i}$ for $1\leq i\leq
N$ and integers $k_{i}\geq1$. Hence, we can associate a degree to monomials
$e_{i_{1}}\dots e_{i_{l}}$ with $(i_{1},\dots,i_{l})\in\{1,\dots,N\}^{l}$ for
$l\geq0$ (note that the monomial associated to the empty sequence is by
convenience $1$), namely
\[
\deg(e_{i_{1}}\cdots e_{i_{l}})=\sum_{r=1}^{l}k_{i_{r}}.
\]
We denote by $\mathbb{A}_{\deg\leq m}^{N}$ the vector space of polynomials
generated by monomials of degree less or equal $m$, for some integer $m\geq1$.
We define a continuous map $\phi:\mathbb{R}^{N}\rightarrow\mathbb{A}_{\deg\leq
m}^{N}$, via
\[
\phi(x_{1},\dots,x_{N})=\sum_{l\geq0}\sum_{_{(i_{1},\dots,i_{l})\in
\{1,\dots,N\}^{l},\ \sum_{r=1}^{l}k_{i_{r}}\leq m}\ }x_{i_{1}}\cdots x_{i_{l}%
}e_{i_{1}}\cdots e_{i_{l}}.
\]
Continuity is obvious, since we are given monomials in each coordinate. $\phi$
is even an embedding and a closed map.

The following example shows the relevant idea in coordinates, since for $N=1$
and $\deg(e_{1})=1$ we obtain $\mathbb{A}_{\deg\leq m}^{1}=\mathbb{R}^{m+1}$.

\begin{example}
\label{example:1}Fix $m\geq1$. Then $\phi(x)=(1,x,x^{2},\dots,x^{m})$ is a
continuous map $\phi:\mathbb{R}^{1}\rightarrow\mathbb{R}^{m+1}$. Given a
positive Borel measure $\mu$ on $\mathbb{R}^{1}$ such that moments up to
degree $m$ exist, i.e.
\[
\int_{\mathbb{R}}|x|^{k}\mu(dx)<\infty
\]
for $0\leq k\leq m$, then $\phi_{\ast}\mu$ admits moments up to degree $1$.
Hence we conclude that there exist $1\leq k\leq m+1$, points $x_{1}%
,\dots,x_{k}$ and weights $\lambda_{1},\dots,\lambda_{k}>0$ such that
\[
\int_{\mathbb{R}^{N}}P(x)\mu(dx)=\sum_{i=1}^{k}\lambda_{i}P(x_{i})
\]
for all polynomials $P$ of degree less or equal $m$.
\end{example}

\begin{theorem}
\label{full cubature}Given $N\geq1$ and degree function $\deg$ and $m\geq1$.
Fix a finite, positive Borel measure $\mu$ on $\mathbb{R}^{N}$ concentrated in
$A\subset\mathbb{R}^{N}$, i.~e.~$\mu(\mathbb{R}^{N}\setminus A)=0$, such that
\[
\int_{\mathbb{R}^{N}}\abs{x_{i_{1}}\cdots x_{i_{l}}}\mu(dx)<\infty
\]
for $(i_{1},\dots,i_{l})\in\{1,\dots,N\}^{l}$ with $\sum_{r=1}^{l}k_{i_{r}%
}\leq m$. Then there exist an integer $1\leq k\leq\dim\mathbb{A}_{\deg\leq
m}^{N}$, points $x_{1},\dots,x_{k}\in A$ and weights $\lambda_{1}%
,\dots,\lambda_{k}>0$ such that
\[
\int_{\mathbb{R}^{N}}P(x)\mu(dx)=\sum_{i=1}^{k}\lambda_{i}P(x_{i})
\]
for $P\in\mathbb{A}_{\deg\leq m}^{N}$.
\end{theorem}

\begin{proof}
The measure $\phi_{\ast}\mu$ admits first moments by assumption, hence we
conclude by Corollary \ref{cubature result cooked down}.
\end{proof}

\begin{remark}
Tchakaloff's Theorem is a special case of the above theorem with $A = \supp
\mu$.
\end{remark}

\begin{remark}
Fix a non-empty, closed set $K\subset\mathbb{R}^{N}$. We note that a finite
sequence of real numbers $m_{i_{1}\dots i_{l}}$ for $(i_{1},\dots,i_{l}%
)\in\{1,\dots,N\}^{l}$ with $\sum_{r=1}^{l}k_{i_{r}}\leq m$ represents the
sequence of moments of a Borel probability measure $\mu$ with support
$\supp\mu\subset K$, where moments of degree less or equal $m$ exist, if and
only if
\[
\sum_{l\geq0}\sum_{_{(i_{1},\dots,i_{l})\in\{1,\dots,N\}^{l},\ \sum_{r=1}%
^{l}k_{i_{r}}\leq m}}m_{i_{1}\dots i_{l}}e_{i_{1}}\cdots e_{i_{l}}\in\conv
\phi(K).
\]
The argument in one direction is that any element of $\conv\phi(K)$ is represented as expectation with respect to
some probability measure with support in $K$, for instance the given convex combination. The other direction is
Tchakaloff's Theorem in the general form of Theorem~\ref{full cubature}. Consequently we have a precise geometric
characterization of solvability of the Truncated Moment Problem for measures with support in $K$. Notice that one
can often describe $\conv\phi(K)$ by finitely many inequalties.
\end{remark}

\end{document}